\documentclass[11pt]{amsart}
\usepackage{amsmath}
\usepackage{amssymb}
\usepackage{amsthm}
\usepackage{amscd, amsfonts}
\usepackage[dvips]{epsfig}
\usepackage{verbatim}

\usepackage{epsf}
\usepackage[usenames]{color}
\usepackage[bookmarksnumbered,pdfpagelabels=true,plainpages=false,colorlinks=true,
            linkcolor=black,citecolor=black,urlcolor=blue]{hyperref}

\newtheorem{theorem}{Theorem}[section]
\newtheorem{lemma}[theorem]{Lemma}

\theoremstyle{definition}

\newtheorem{remark}[theorem]{Remark}
\newtheorem{convention}[theorem]{Convention}

\numberwithin{equation}{section}

\newcommand{\Si}{\Sigma}
\newcommand{\al}{\alpha}
\newcommand{\be}{\beta}
\newcommand{\cF}{\mathcal{F}}
\newcommand{\rc}{\mathrm{Ricc}}

\begin{document}

\title[Einstein-Euler]{On the Choquet-Bruhat--York--Friedrich
 formulation of the Einstein-Euler equations}

\author[Disconzi]{Marcelo M. Disconzi}
\address{Department of Mathematics,
Vanderbilt University,  1326 Stevenson Center, Nashville, TN 37240, USA.}
\email{marcelo.disconzi@vanderbilt.edu}
\thanks{Marcelo M. Disconzi is supported by NSF grant 1305705.}

\author[Pingali]{Vamsi P. Pingali}
\address{Department of Mathematics,
Johns Hopkins University, 404 Krieger Hall
 3400 N. Charles Street, Baltimore, MD 21218, USA.}
\email{vpingali@math.jhu.edu}



\begin{abstract}
Short-time existence for the Einstein-Euler and the vacuum Einstein equations is proven using a Friedrich inspired formulation due to Choquet-Bruhat and York,
where the system is cast into a symmetric hyperbolic form and the Riemann
tensor is treated as one of the fundamental unknowns of the problem. The reduced system of Choquet-Bruhat and York, along with the preservation of the gauge, is shown to imply the full Einstein equations.
\end{abstract}

\maketitle

\section{Introduction}
In the vast amount of literature that exists on the Cauchy problem of General
Relativity (GR)\footnote{A complete or extensive account
of all references is beyond the scope of this paper, whose
length we tried to keep short. We refer the interested reader to the monographs
\cite{CBbook, Ringstrom} and the survey papers \cite{CGP, FriRen}. A long, although far
from complete, review of the literature of the Cauchy problem for the Einstein-Euler
system specifically, is given in \cite{Dis}, while a thorough and up-to-date treatment of relativistic fluids can be found in \cite{RZ}. Further discussion
on relativistic fluids, including problems such as the inclusion of  viscosity,
long-time existence, and other fluid-matter models, can be found in 
\cite{Dis2, RodSpeck, Speck1, Speck2}
and references therein.},
the formulation in terms of a
first order symmetric hyperbolic system (FOSH) has recently attracted significant
attention (see e.g. \cite{CBY, FriHypRed, F, FriNagy, FriRen} and references therein.).
Here, we focus on the Choquet-Bruhat and York \cite{CBY} formulation of the Einstein-Euler system in terms of the
Lagrangian\footnote{Intuitively, the picture is like this.
We can think of two ways to study the flow of a river: one could float downstream on a boat, or one could sit on the bank and observe the flow. The former (Lagrangian description) corresponds to tracking the position of every particle and the latter (Eulerian description) to observing the velocity vector field. Both descriptions are useful in the study of relativistic and non-relativistic fluids.}
description
of the fluid flow, which itself was adapted from an earlier formulation by Friedrich \cite{F}.

In \cite{CBY} the authors wrote a system of equations in terms of the Riemann tensor (as opposed to the Weyl tensor used in \cite{F}), and chose a gauge that reduced this system to a FOSH. It turned out that this system had physical characteristics (in contrast with Friedrich's one), i.e. the assumption that the speed of sound in the fluid is less than that of light was crucial to prove the hyperbolicity of the equations. This is important, for example,
because it gives a natural breakdown criterion for the problem of long-time existence. However, the task of proving that the gauge is preserved and the original Einstein-Euler system is satisfied was not carried out in \cite{CBY}. In this article, we complete the proof of short-time existence for the Einstein-Euler system \`a la Choquet-Bruhat and York \cite{CBY}. In what follows, we shall
restrict ourselves to barotropic fluids.

In mentioning Choquet-Bruhat and York's construction of a FOSH, it is worth
recalling the general strategy for solving Einstein's equations.
The Einstein equations do not form an ``honest" system of evolution equations,
in the sense that some of the equations are constraint equations.
Such a difficulty is a consequence of the diffeomorphism invariance,
or gauge freedom, enjoyed by the system.
To circumvent this problem one considers a different set of equations,
generally referred to as the ``reduced system" containing
 only genuine evolution equations that can be solved
using standard techniques. This system is chosen so as to correspond
to the original Einstein's equations \emph{modulo} the constraints 
(which have to be solved separately in order to produce
a full set of initial data for the evolution problem, see below).
This task can be accomplished by a suitable choice of gauge.
A solution to the original system is then obtained by
showing that the gauge conditions are in fact satisfied
on the time-interval where a solution to the
reduced equations has been shown to exist, \emph{provided}
they are satisfied initially, i.e., at time $t=0$. This is done
by deriving a suitable system of evolution equations for the
gauge and using uniqueness.

We remark that the result here obtained is not, in itself, new.
Short-time existence for the Einstein-Euler
system had been proven earlier by Choquet-Bruhat \cite{C2}, and subsequently
by Lichnerowicz \cite{Lich_fluid_1,Lich_fluid_2}. The novelty in the approach
initiated by Friedrich \cite{F} is the use of the Lagrangian
description of the fluid.
This sheds new light in the problem of the so-called  ``fluid body" modeling
certain stellar dynamics, where one
attempts to solve the free-boundary problem that arises from considering the
the system formed by the Einstein equations coupled to the Euler equations within a bounded
region, with  vacuum Einstein's equations holding on the complement.
Recent existence results for this problem have been obtained by
Brauer and Karp \cite{BrauerKarp1, BrauerKarp2}.

\section{Summary of results}
\indent In the study of the Cauchy problem in GR, one is usually given a Riemannian smooth $3$-fold $(\Sigma, h_0)$, a symmetric $2$-tensor $K$, and other initial data corresponding to the matter fields. This initial data is required to satisfy certain constraint equations, which are derived from the
Gauss-Codazzi-Mainardi equations and the Einstein equations,
and ensure that $(\Sigma, h_0)$ embeds isometrically, with $K$ as
its second fundamental form, into 
the space-time that is eventually obtained as a solution of the full
Einstein system.
 This prescription of data is usually facilitated by means of the conformal method of solving the constraints. The aim then is to find an Einsteinian development, i.e. a Lorentzian $4$-fold $(M,g)=(\Sigma \times [0,T] , g)$, containing matter fields satisfying Einstein's equations and obeying the initial conditions on the matter fields.

  Naturally, upon writing $\Sigma \times [0,T]$ we are relying on a particular choice
 of
diffeomorphism to parametrize the would-be\footnote{As there is no natural notion
of a time coordinate in GR.} ``time coordinate" $t \in [0,T]$. Although the existence
of a solution to the Einstein-Euler system can be stated in a more invariant fashion,
here it is convenient to write explicitly $\Si \times [0,T]$ in order to
 follow the similar statements of \cite{CBY}, on which this work is largely based,
 and also to facilitate the identification of the spaces where solutions live in.

From a PDE perspective,
 using a standard 3+1 coordinate decomposition where the vectors
 $\frac{\partial}{\partial x^i}$, $i=1,2,3$ are space-like and
 $\frac{\partial}{\partial x^0}$ is time-like, the 
constraint equations read  \cite{CBbook}
\begin{gather}
 \mathrm{Ricc}_{\mu 0} - \frac{1}{2} R g_{\mu 0} =  T_{\mu 0}
  \text{ on } \Si,
\label{constraints}
\end{gather}
where $T$ is the  stress-energy tensor of the matter fields.
It is easy to see that
 the constraints do not form a system of second order 
evolution equations. In particular, initial data for the full Einstein system 
ought to satisfy (\ref{constraints}), and thus cannot be given 
arbitrarily. From these considerations, it is seen that, 
while the construction of initial data 
satisfying the constraint equations is doubtless a crucial aspect of 
the investigations surrounding  Einstein's equations,
it can be considered apart from the  evolution problem.
Thus, in what follows, it is assumed throughout that
 the fields in a given initial data set \emph{always satisfy the constraint equations}.
 We comment further on the
 initial conditions in section \ref{initial_data_section}.

We assume that we are given the aforementioned type of initial data. This data is then converted into the type we need for solving the FOSH. In what follows, for the fluid case, $p\geq 0$ indicates pressure, $\mu(p) > 0$ is the density as a function of $p$, $v$ is the initial $3$-velocity of the fluid on $\Sigma$, and $u$ is the $4$-velocity field on $M$. We prove our results in the Sobolev spaces $H^s$.
In what follows,  $\Si$, $h_0$, and $p$ denote the quantities just described, and
repeated indices are summed over. We also assume the reader is familiar
 with the terminology and the Cauchy problem in  GR and the initial
 conditions for the Einstein-Euler system.

Our first result is on the vacuum Einstein equations, i.e. $\mathrm{Ricc}=0$:
\begin{theorem}
Let  $s>\frac{3}{2}+2$, and let $(\Sigma, h_0, K)$ be an initial data set for the vacuum
Einstein equations, with $h_0$ in
$H^{s+1} (\Sigma)$, $K_0$ in $H^{s} (\Sigma)$, and $\Si$ compact. Then,
there exists an Einsteinian development $M=(\Si \times [0,T],g)$ satisfying $\mathrm{Ricc}(g)=0$. The metric $g$ thus obtained is in $C^{0}([0,T],H^{s+1}(\Si))\cap C^{1}([0,T],H^{s}(\Si)) \cap C^{2}([0,T],H^{s-1}(\Si))$.
\label{vac}
\end{theorem}

We have assumed  $\Si$ to be compact for simplicity. This can be relaxed provided suitable
asymptotic conditions on the initial data are imposed. While asymptotically flat initial data
is the standard choice for the vacuum case, existence under similar conditions
becomes technically challenging for the case of the Einstein-Euler system;
see \cite{BrauerKarp1,BrauerKarp2}. Notice, also, that by stating our existence theorem on the
\emph{closed} interval $[0,T]$, we are not taking the maximal Cauchy development of the initial data.

For perfect fluids, the energy-momentum tensor is $\mathbf{T} = (\mu +p) \mathbf{u}\otimes \mathbf{u} + p \mathbf{g}$. The Einstein-Euler system for a perfect fluid is
\begin{align}
\begin{split}
& \mathrm{Ricc}_{\mu \nu} - \frac{1}{2} R g_{\mu \nu} =  T_{\mu \nu}  \\
& (\mu+p) u^{\al} \nabla _{\al} u ^{\be} + (u^{\al} u^{\be}+g^{\al \be})\partial _{\al} p = 0  \\
& (\mu+p) \nabla _{\al} u ^{\al} + u^{\al}\partial _{\al} \mu = 0  \\
& u^{\al} u_{\al} = -1
\end{split}
\label{EinEul}
\end{align}
Notice that we have chosen units such that $\frac{8\pi G}{c^2} = 1$.
Also note that the first equation maybe written as
$\mathrm{Ricc} _{\mu \nu} = \rho _{\mu \nu}$ where $\rho _{\mu \nu} = T_{\mu \nu} - \frac{\mathrm{Tr}(\mathbf{T})}{2} g_{\mu \nu}$, where $\mathrm{Tr}$ denotes the trace.

The Einstein-Euler system also exists for a short period of time:
\begin{theorem}
Let $s>\frac{3}{2}+2$, and let $(\Sigma, h_0, K, p_0, v)$ be an initial data set
for the Einstein-Euler system,
 where $\Sigma$ is compact, $h_0$ is in $H^{s+1}(\Sigma)$, $K_0$ is in $H^{s}(\Sigma)$, $p_0$ is in $H^{s}(\Sigma) $, and $v$ is in
 $H^{s} (\Sigma)$. Fix a smooth invertible function $\mu :[0,\infty)\rightarrow (0,\infty)$
with $\mu^{\prime} \geq 1$. Then there exists an Einsteinian development $M=(\Si  \times [0,T],g)$ satisfying the Einstein-Euler system. The metric $g$ thus obtained is in $C^{0}([0,T],H^{s+1}(\Si ))\cap C^{1}([0,T],H^{s}(\Si)) \cap C^{2}([0,T],H^{s-1}(\Si ))$, the four-velocity $u \in C^{0}([0,T],H^{s}(\Si))\cap C^{1}([0,T],H^{s-1}(\Si ))$, and the pressure $p\in C^{0}([0,T],H^{s}(\Si ))\cap C^{1}([0,T],H^{s-1}(\Si )) $. They obey the initial conditions, hence in particular
 the orthogonal projection of $u$ onto $T\Sigma$ is $v$.
\label{EE}
\end{theorem}

\begin{remark}
The condition $\mu^{\prime} \geq 1$ guarantees that the speed of sound is at most that
of light. The Einstein-Euler system will be a FOSH only as long as this condition holds.
This is one advantage of having only physical characteristics, as mentioned in the introduction.
\end{remark}

\begin{remark}
Following the usual arguments relying on the finite-propagation speed property of FOSH systems,
in the proof of theorems \ref{vac} and \ref{EE} we shall work solely on a single coordinate chart and use uniqueness of solutions of FOSH systems in the following way: Let $V \subset \tilde{U}  \subset U \subset \Sigma$ be relatively compact open sets having smooth boundary (with $V$ and $\tilde{U}$ contained in a coordinate chart). The domain of a solution to a FOSH system with initial data in $\tilde{U}$ contains $V \times [0,T_{V}]$ for some $T_{V}$ such that $V \times [0,T_{V}]$ in   the domain of influence of $\tilde{U}$. Moreover, if  $V \cap W$ is not empty for a relatively compact, open, smooth $W$ then by uniqueness of solutions the solutions on $V \times [0,\min{T_V, T_W}]$ and $W \times [0,\min{T_V, T_W}]$ coincide. This way we get a unique solution on $U\times [0,T_U]$ for some $T_U>0$. A solution on the whole of $\Si \times [0,T]$, for some
$T>0$, is then obtained by a standard gluing procedure.
\label{remark_charts}
\end{remark}

\begin{remark}
The regularity hypotheses in the theorems along with Sobolev embedding imply that the metric is $C^3$. This would appear to be superfluous because Einstein's equations involve only two derivatives of the metric.
 However, the reduced system involves the derivatives of the curvature. Hence the need for additional regularity. In fact, in whatever follows, we shall need to take two derivatives of the curvature.
 We are so allowed because $h_0$ belongs to $H^{s+1}$ and thus has four (weak) derivatives since
 $s > \frac{3}{2} + 2$.
\end{remark}

 As already mentioned, Einstein's equations (in vacuum and coupled to matter)
 are diffeomorphism invariant and hence an appropriate gauge has to be chosen in order to solve them. Traditionally, harmonic coordinates were employed to convert the Einstein equations into a second order hyperbolic system. However, in the Lagrangian framework, Choquet-Bruhat and York chose the so-called Cattaneo-Ferrarese (CF) gauge consisting of Lagrangian observers following the fluid flow. In other words, a choice of a local orthonormal frame $\{ e_{\al} \}_{\al=0}^3$ such
 that $e_0 =  u = \partial _{x^0}$, and the remaining $\{ e_i \}_{i=1}^3$ are
 Fermi propagated, i.e., $\omega ^{j} _{0i}=0$ for $i,j=1,2,3$, where $\omega$ denotes
 the connection coefficients of the Levi-Civita connection (or Ricci rotation coefficients)
  with respect to
  $\{ e_{\al} \}_{\al=0}^3$. It is also assumed that
local coordinates $x^{\al}$ have been chosen such that
  $e_0 \equiv u = \partial _{x^0}$ and $\partial_{x^i}$ gives a basis for the tangent
  space of $\Si$ within the selected coordinate chart. For notational convenience, we
  shall denote $f_{ij} = \omega ^{j}_{0i}$ and  $\partial _i = e_i$, $i=1,2,3$.
From the above, we can write
\begin{gather}
 \partial _i \equiv e_i = A_{i}^j \left (\partial _{x^j}-b_j \partial _{x^0} \right ),
 \nonumber
 \end{gather}
 for $i=1,2,3$, and a certain invertible matrix $A$ and a one form $b$. Note that $\partial _{\al}$ should not be confused with a coordinate basis
 (which are denoted by $\partial _{x^{\al}}$).
Notice, also, that $\omega$ satisfies (see also lemma \ref{initsymm})
\begin{gather}
\omega_{ij}^k = -\omega_{ik}^j \text{ and }
\omega_{ij}^0 = \omega_{i0}^j,
\nonumber
\end{gather}
for $i,j,k=1,2,3$.  In this gauge, the Einstein equations were re-written to form a reduced FOSH \cite{CBY}.

\begin{convention}
From now on, Latin indices run from 1 to 3 and Greek indices from 0 to 3.
\end{convention}

Following \cite{CBY}, all the symbols appearing henceforth are to be treated as ``abstract" --- for example, $R_{\al\be \mu \nu}$ is not known, a priori, to be the Riemann tensor of a metric; except, however, for those quantities determined at $t=0$, in which
case they do have their usual meaning. The strategy is to write evolution equations for these ``abstract" quantities, identify them as FOSH systems, and use uniqueness to conclude that indeed these ``abstract" symbols correspond to the ``correct" geometric objects.
For the sake of brevity, we use the symbol $\nabla$ or a semicolon to mean ``covariant derivative'' with the ``correct'' connection coefficients in the chosen frame i.e. $\omega ^{a} _{0b} =0, \omega ^{0}_{0a} = Y_a = \omega ^{a}_{00}, \omega ^{0} _{ab} = X_{ab} = \omega ^{b} _{a0}$. We denote spatial covariant derivatives and curvatures with a $\widetilde{\omega}$ and $\widetilde{R}$ respectively. Square brackets enclosing two letters $A_{[a,b]} = A_{ab} - A_{ba}$ indicates antisymmetrisation, except
that we leave out a conventional factor of $\frac{1}{2}$,  whereas angular ones enclosing three letters separated by commas ($A_{<a,b,c>}$) indicates cyclic summation $A_{abc}+A_{bca}+A_{cab}$.
This notation is not standard but is useful in this context.

In the case of vacuum, there are no fluid flow lines. Hence we may impose the additional gauge choice $Y_i=0$. The reduced system of equations as written in \cite{CBY} is
\begin{align}
\begin{split}
& \partial _0 a^{i} _{j} = -a^k _j X _{ik}   \\
& \partial _0 b_i = 0   \\
& R^{\_ \_ i}_{0h\_j} = \partial _0 \omega ^i _{hj} + X_{h}^k \omega ^{i} _{kj}  \\
& R_{h0i0} = -\partial_0 X_{hi} - X_{h}^j X_{ji}  \\
& \nabla _0  R_{hk \lambda \mu} = -\nabla _k R_{0h \lambda \mu} +\nabla _h R_{0k\lambda \mu}  \\
& \nabla _0 R_{0h \lambda \mu} = \nabla ^{l} R_{lh \lambda \mu}
\end{split}
\label{reducedvac}
\end{align}
where $a= A^{-1}$. In the above, and in what follows, we adopt the following
notation: underbars ``$\_$" are used to denote empty slots in the order
of the indices
when one raises or lowers an index. For instance, in $R^{\_ \_ i}_{0h\_j}$
the two first $\_$'s on the top and the $\_$ on the third entry on the bottom
tell us that the upper index $i$ was obtained by raising the third lower index
from $R_{0hij}$. Although this notation is not completely standard, it is similar
to the one used in \cite{CBY}, which we tried to follow.

In the perfect fluid case, let $F = \int \frac{dp}{\mu (p) +p}$, $\rho _{00} = \frac{1}{2} (3p + \mu)$,$\rho _{i0}=0$, and $\rho_{ij} = \delta _{ij} \frac{1}{2}(\mu-p) $. The reduced system is
\begin{align}
\begin{split}
&\partial _0 a^{i} _{j} = -a^k _j X _{ik}   \\
&\partial _0 b_i = -a_i ^h Y_h   \\
 &\partial _0 \omega ^i _{hj} + X_{h}^k \omega ^{i} _{kj} + Y^i X_{hj} -Y_j X_{h\_}^{\_ i} = R^{\_ \_ i}_{0h \_ j}  \\
 & \partial_0 X_{hi} + X_{h\_}^{\_ j} X_{ji}-Y_h Y_i-\widetilde{\nabla}_i Y_h - \frac{X_{l}^l}{\mu ^{\prime}}(X_{hi}-X_{ih}) = -R_{h0i0} \\
&\mu ^{\prime} \partial _0 Y_h - \widetilde{\nabla} ^l X_{hl}-Y^l(X_{li}-X_{il}) + \mu ^{\prime} (Y_h \partial _0 F -X_{h\_}^{\_ l}\partial _l F) \\
& \, +\partial _h \mu ^{\prime} \partial _0 F = 0  \\
&\nabla _0  R_{hk \lambda \mu} = -\nabla _k R_{0h \lambda \mu} +\nabla _h R_{0k\lambda \mu}  \\
&\nabla _0 R_{0h \lambda \mu} = \nabla ^l R_{lh \lambda \mu} + \nabla _{\mu} \rho _{\lambda h}-\nabla _{\lambda} \rho _{\mu h}  \\
& \partial _0 \mu = -(\mu+p) X_{l}^l
\end{split}
\label{reducedEE}
\end{align}
In both vacuum and fluid cases, the constraints are obtained 
from the splitting of the Riemann tensor of $(M,g)$ into
that of $\Si$ and the second fundamental form of $\Si$ inside $M$
(in other words, from the Gauss-Codazzi-Mainardi equations)
and Einstein's equations, upon restriction to 
$\Si = \{ t = 0\}$. This enables us to solve for temporal derivatives of all the quantities at $t=0$. Substituting these expressions in the so-called quasi-constraints\footnote{This procedure is necessary because $\partial _i$ contains $\partial _{x^0}$.} which we write below, gives us the actual constraints 
(see also \cite{CBY}). 
The initial connection is the Levi-Civita one. The other quasi-constraints are
\begin{align}
\begin{split}
&\nabla _{\langle h,} R_{i,j\rangle \lambda \mu} = 0  \\
&\nabla _h R_{h0\lambda \mu} = \nabla _{\lambda} \rho _{\mu 0} - \nabla _{\mu} \rho _{\lambda 0}  \\
&R_{hkij} = \tilde{R} _{hkij} + X_{ki}X_{jh}-X_{jk}X_{ih}  \\
&-R_{kh0j} = \widetilde{\nabla} _{[k,}X_{h]j}-Y_j X_{[k,h]}
\end{split}
\label{qconsEE}
\end{align}

\section{Initial data\label{initial_data_section}}
The initial data required for the reduced system is defined on $\Sigma \times \{0\}$:
\begin{itemize}
\item A field of coframes $a_{j} ^i$ and of covectors $b_i$ creating a metric $h_0$ on $\Sigma$ via $h_0 ^{jh}=a^j _{l} a ^{lh} -b^j b^h$ which is assumed to be positive definite. The initial metric on the manifold $M$ is $g (t=0)=-(\theta ^0) ^2 + \sum (\theta ^i) ^2$ where $\theta ^i = a^i _{j} dx^j$ and $\theta ^0 = dx^0 + b_i dx^i$.
\item Fields $\omega _{ij} ^k$, $X_{ij}$, and $Y_i$ (which is assumed to be zero in the vacuum case). These are supposed to define the connection coefficients of the Levi-Civita connection of $g$ initially (with $f_{ij}=0$).
\item Tensor components $R_{ijkl}$, $R_{0i jl}$, and $R_{0i0j}$ that define the Riemann curvature tensor initially.
\item In the case of the perfect fluid, we also need $\mu(p)> 0$ obeying $\mu^{\prime}\geq 1$, and $p\geq 0$.
\end{itemize}
In addition, the Einstein equations are imposed on this initial data at $t=0$
in order to derive the relation between the usual Eulerian initial data of the $3+1$ decomposition --- which is given  in theorems \ref{EE} and \ref{vac} ---
and the initial data needed for the FOSH systems, as we explain below. A detailed
account of the correspondence between initial data sets for the Einstein-Euler
system and those of reduced equation in Lagrangian coordinates can be obtained
by an argument similar to that of \cite{Dis}.

 Given a Riemannian 3-fold $(\Sigma, h_0)$, choose local coordinates $\tilde{x}^i$ on it. Embed it into $M=\Sigma \times \mathbb{R}$ as $\Sigma \times \{0\}$. Let $\tilde{e}_i$ be an orthonormal frame on $\Sigma$ and let $\tilde{\theta}^i$ be the dual coframe. Then $h_0 = \sum \tilde{\theta ^i} \otimes \tilde{\theta ^i}$. In the vacuum case, we may simply define the initial Lorentz metric on $\Sigma \times \{0\}$ as $g=-(dx^0)^2+h_0$. This corresponds to $b_i$ being zero initially.

In the case of a perfect fluid, define a metric $g=h_0 + v_i dx^0 \tilde{\theta ^i} -(dx^0)^2$ on $TM$ restricted to $\Sigma \times \{0\}$ with $v_i$ being the components of the dual (with respect to $h_0$) of $v$. This is a Lorentzian metric, with $e_0=\partial _{x^0}$ being a unit timelike vector projecting
to $v$, and restricting to $h_0$ on $\Sigma \times \{0\}$. Complete $e_0$ to an orthonormal basis $e_{\al}$. This gives us $a_{j} ^i$ and $b_i$ lying in $H^s (\Sigma)$. Calculations similar to the ones in \cite{Dis} maybe used to define the remaining fields, such as $\omega, X$, on $\Sigma \times \{0\}$.

The above reasoning combined with the fact that both (reduced) systems
(\ref{reducedvac}) and
(\ref{reducedEE}) are quasilinear FOSH having initial data in (at least) $H^{s-1}$ implies that both systems have solutions in $C^{0}([0,T_1],H^{s-1}(U))\cap C^{1}([0,T_1],H^{s-2}(U))$,
where $U$ is some local chart as described in remark \ref{remark_charts}. The mismatch between the regularity of the initial data and that of the solution is then corrected by a bootstrap argument as in \cite{Dis} following the results of \cite{FM}.

Thus we have a solution to both systems with $a ^{i} _{j} \in C^0 ([0,T_1], H^{s+1}(U))\cap C^1 ([0,T_1], H^{s}(U)) \cap C^2 ([0,T_1], H^{s-1}(U))$, $b _i$, $\omega ^{k} _{ij}$, $X_{ij}$, $Y_i$, $p$ in $C^0 ([0,T_1], H^{s}(U)) \cap C^1 ([0,T_1], H^{s-1}(U))$, and $R_{\al \be \mu \nu}$ in $C^0 ([0,T_1], H^{s-1}(U)) \cap C^1 ([0,T_1], H^{s-2}(U))$.

\section{Proofs}
 We prove that the constraints and the gauge are preserved. This implies that if the Einstein equations are satisfied initially, then they are satisfied in the future. We accomplish these steps by proving that the relevant quantities satisfy FOSH systems with zero as their unique solution. Note that by definition, $\partial _0 = \frac{\partial}{\partial x^0}$ and $\partial _i = A_i ^{l}(\frac{\partial}{\partial x^l} - b_l \frac{\partial}{\partial x^0})$. We also note that if a linear symmetry of the $R_{\al \be \mu \nu}$ is satisfied initially, then its spatial derivatives are zero. Since the temporal derivatives are related to the spatial ones by the evolution equations (which are imposed on the variables at $t=0$), we see that $\partial _i$ applied to such a symmetry also yields zero.

\subsection{Vacuum}
   Firstly, we see as to why the vacuum Einstein equations are implied by the preservation of the constraints and the gauge:
\begin{lemma}
If the gauge and the constraints are preserved, then $\mathrm{Ricc}_{\alpha \beta}$ is zero in the future, if so initially.
\label{separation}
\end{lemma}
\begin{proof}
Notice that $d\mathrm{scal} = 2 \mathrm{div} \mathrm{Ricc}$, where $\mathrm{scal}$
is the scalar curvature and $\mathrm{div}$ means divergence. The given system implies that (assuming the constraints and the gauge) $\mathrm{Ricc}_{\alpha \beta ; \gamma} = \mathrm{Ricc}_{\alpha \gamma ;\beta} $. Contracting $\alpha$ and $\beta$, we see that $d\mathrm{scal} = 0$ i.e. $\mathrm{scal}=0$ because it is so, initially. Hence
\begin{align}
\begin{split}
\rc_{0i;0} &= \, \rc_{00;i}  \\
\rc_{00;0} &= \, \displaystyle \sum _{i=1} ^{3} \rc_{0i;i}
\end{split}
\label{one}
\end{align}
The leading matrix ($M_0$) for the equations above is
\[ \left( \begin{array}{cccc}
1 & -B_1 & -B_2 & -B_3 \\
-B_1 & 1 & 0 & 0 \\
-B_2 & 0 & 1 & 0 \\
-B_3 & 0 & 0 & 1
\end{array} \right)\]
where $B_i = -A^{j} _i b_j$. It is positive definite (see lemma 11 in \cite{CBY}).
\begin{align}
\begin{split}
\rc_{ij;0} &= \, \rc_{i0;j} \\
\rc_{i0;0} &= \, \rc_{ij}^{;j}
\end{split}
\label{two}
\end{align}
Equations (\ref{one}) and (\ref{two}) form a FOSH system (the leading matrix of equation \ref{two} also has positive eigenvalues by a similar argument as for
equation (\ref{one})). Hence $\rc_{\alpha \beta} = 0$. Note that we treated $\rc_{\alpha \beta}$ and $\rc_{\beta \alpha}$ as distinct variables.
\end{proof}

Now, we prove that the constraints and the gauge are preserved i.e., among other things $\partial _{\al}$ forms an orthonormal frame for a metric whose Levi-Civita connection's components are $\omega ^{0}_{ab} =X_{ab}$, $\omega _{0\al} ^{\be} =0$, $\omega _{ij} ^k$ and whose Riemann curvature tensor is $R_{\al \be \mu \nu}$. For further use, we prove some symmetries of $R_{\al \be \mu \nu}$:
\begin{lemma}
The following relations are satisfied for some time:
\begin{eqnarray}
R_{\al \be \mu \nu} &=& -R_{\al \be \nu \mu} \nonumber \\
R_{\al \be \mu \nu} &=& -R_{\be \al \mu \nu} \nonumber \\
\omega ^{p} _{ij} &=& -\omega ^{j} _{ip} \nonumber
\end{eqnarray}
\label{initsymm}
\end{lemma}
\begin{proof}
\begin{eqnarray}
\partial _{0} (R_{hk0j}+R_{kh0j}) &=& -\nabla_k R_{0h0j}+\nabla_{h}R_{0k0j}+\nabla_k R_{0h0j}-\nabla_{h}R_{0k0j} \nonumber\\
&=& 0 \nonumber
\end{eqnarray}
Hence $R_{hk0j}=-R_{kh0j}$ (since it is so, initially). Similarly, $R_{ijkl}=-R_{jikl}$. We also have
\begin{eqnarray}
\partial _{0} (R_{ijkl}+R_{ijlk}) &=& -\nabla_j(R_{0ikl}+R_{0ilk})+\nabla_i(R_{0jkl}+R_{0jlk}) \nonumber \\
\partial _{0} (R_{0ikl}+R_{0ilk}) &=& \nabla ^p (R_{pikl}+R_{pilk}) \nonumber
\end{eqnarray}
The system above is FOSH (having zero as its unique solution). Indeed, the leading matrix is
\[ \left( \begin{array}{cccccc}
1 & 0 & 0 & B_2 & -B_1 & 0 \\
0 & 1 & 0 & B_3 &0& -B_1 \\
0 & 0 & 1 & 0 & B_3 &-B_2 \\
B_2 & B_3 & 0 & 1 & 0 & 0 \\
-B_1 & 0 &B_3 & 0 & 1 & 0 \\
0 & -B_1 & -B_2 & 0 & 0 & 1
\end{array} \right)\]
Its eigenvalues are $1, 1\pm \sqrt{B_1 ^2 + B_2 ^2 + B_3 ^2}$ with multiplicity $2$. They are positive for some time (by the assumptions on $a$ and $b$). This means that $R_{ijkl}=-R_{ijlk}$ and $R_{0ikl}=-R_{0ilk}$. Using these symmetries of $R_{\alpha \beta \mu \nu}$, we see that $\partial _{0} (\omega _{ij} ^{k} + \omega _{ik} ^j) = X_i ^{l} (\omega _{lj} ^{i}+\omega_{li} ^{j})$. Hence $\omega _{ij} ^{k} = -\omega _{ik} ^{j}$.
\end{proof}
\indent By explicit calculation, we see that $\partial _0 \partial _i - \partial _i \partial _0 = -X_i ^{a} \partial _a$ and that
\begin{eqnarray}
\partial _{[i,}\partial _{j]} &=& (2f_{ij}+X_{[i,j]})\partial_0 + c_{ij} ^p \partial _p \nonumber
\end{eqnarray}
where $f$ and $c$ satisfy $f_{ij}=0$ and $c_{ij}^p=\omega_{[i,j]}^p$ when $t=0$. The evolution equations for $f_{ij}$ and $v_{ij}^p=c_{ij}^p-\omega_{[i,j]}^p$ are obtained by differentiating the above equation and using the main evolution system.
\begin{align}
\partial_0 (v_{ij}^p) \partial_p = & \, \partial _0 \partial _{[i,}\partial_{j]}-c_{ij} ^p \partial _0 \partial _p -2\partial_0 (f_{ij}) \partial _0 -\partial_0 (X_{[i,j]})\partial _0
-(2f_{ij}+X_{[i,j]})\partial _0 ^2 \nonumber \\
&+ R_{0[i,j]p} \partial _p + \omega_{l[j} ^p X_{i]}^l \partial _p \nonumber \\
 = & \, \partial _{[i,} \partial _0 \partial_{j]} -X_{[i,} ^{a} \partial _{a} \partial _{j]} - c_{ij} ^p (\partial _p \partial _0 - X_p ^a \partial _a)
 -2\partial_0 (f_{ij}) \partial _0 \nonumber \\
 & -(2f_{ij}+X_{[i,j]})\partial _0 ^2
 -\partial_0 (X_{[i,j]})\partial _0  + R_{0[i,j]p} \partial _p + \omega_{l[j} ^p X_{i]}^l \partial _p
\nonumber \\
= & \, \partial _{[i,} \partial_{j]} \partial _0 - \partial _{[i,}(X_{j]} ^b \partial _b) - X^a_{[i,} (\partial _{j]}\partial _a + c^p _{aj} \partial _p + (2f_{aj}+X_{[a,j]})\partial _0)
\nonumber \\
&- c_{ij} ^p (\partial _p \partial _0 - X_p ^a \partial _a)
- 2\partial_0 (f_{ij}) \partial _0 -(2f_{ij}+X_{[i,j]})\partial _0 ^2
\nonumber \\
&  - \partial_0 (X_{[i,j]})\partial _0 + R_{0[i,j]p} \partial _p + \omega_{l[j} ^p X_{i]}^l \partial _p
\nonumber \\
= &\,  (-\partial _{[i,}(X_{j]} ^p)-X^a _{[i}c^p _{a ,j]} + c^{a}_{ij}X _{a} ^p + R_{0[i,j]p} + \omega_{l[j} ^p X_{i]}^l) \partial _p \nonumber \\
& + ( - X^{a}_{[i,}(2f_{aj}+X_{[a,j]}) - 2\partial_0 (f_{ij}) - R_{0[i,j]0} + X^{l} _{[i} X_{l,j]} ) \partial _0 \nonumber
\end{align}
Comparing coefficients we see that
\begin{align}
\begin{split}
\partial_0 (v_{ij}^p) &= \, -\partial _{[i,}(X_{j]} ^p)-X^a _{[i}c^p _{a ,j]} + c^{a}_{ij}X _{a} ^p + R_{0[i,j]p} + \omega_{l[j} ^p X_{i]}^l   \\
\partial_0 (2f_{ij}) &=\, - X^a_{[i,}(2f_{aj}+X_{[a,j]}) - R_{0[i,j]0} + X^l _{[i} X_{l,j]}
\end{split}
\label{evvf}
\end{align}
The system (\ref{evvf}) is a FOSH system for $v$ and $f$ and hence has a unique solution. If the first Bianchi identity and the defining equation of the Riemann tensor are satisfied, then $f=v=0$ is a solution, and hence as promised, $X$ and $\omega$ form the Levi-Civita connection of the metric defined by $\partial _{\al}$. \\
\indent Next, we write evolution equations for the first Bianchi identity. In what follows $R_{\langle \al,\be,\mu \rangle\nu} = R_{\al\be\mu \nu} + R_{\be\mu \al\nu} + R_{\mu \al\be \nu}$.
\begin{align}
\begin{split}
\partial _0 R_{\langle h,c,d \rangle\al} &= \, -\nabla _h R_{\langle c,0,d \rangle\al} - \nabla _d R_{\langle h,0,c \rangle\al} - \nabla _c R_{\langle 0,h,d \rangle\al} + \nabla _{\langle h,}R_{d,c \rangle 0\al}  \\
\partial _0 R_{\langle 0,h,c \rangle\alpha} &= \, \nabla ^l R_{\langle l,h,c \rangle\al} - \nabla ^l R_{hcl\al} + \nabla _0
 R_{hco\al}
 \end{split}
\label{bian1}
\end{align}
The above system is FOSH for the variables $R_{\langle \al, \be, \mu \rangle\nu}$. Indeed, it is symmetric and the eigenvalues of the leading matrix are $1$ and $1\pm \sqrt{B_1 ^2 + B_2 ^2 + B_3 ^2}$ (which are positive). We will write evolution equations for the other terms in the system (\ref{bian1}). This will prove that the unique solution to the above system is $0$ (since it is $0$ initially). \\
\indent Now, we write the evolution equations for  $\nabla ^b R_{hcb\al} -\nabla _0 R_{hc0\al}$ and prove that zero is their only solution. To accomplish this, we ought to prove that the lower order terms in these equations vanish assuming that all the identities (including the Bianchi identities, the Einstein equations, $\nabla ^b R_{hcb\al} -\nabla _0 R_{hc0\al}=0$, $f_{ij}=0$, $c_{ij} ^p = \omega ^p _{ij}-\omega ^p _{ji}$, etc) hold to order zero.
\begin{align}
\begin{split}
& \nabla _0 (\nabla ^b R_{hcb0}) =  \\
& \nabla _0 (\partial ^b R_{hcb0} - \omega _{bh} ^{\al} R^{\_ \_ b}_{\al c \_ 0} - \omega _{bc} ^{\al} R^{\_ \_ b}_{h \al \_ 0} - \omega _{bb} ^{\al} R_{hc\al 0} - \omega _{b0} ^{\al} R^{\_ \_ b}_{hc\_\al} )   \\
=& \, [\partial _0, \partial ^b] R_{hcb0} - \partial ^b \nabla _c R_{0hb0} + \partial ^b \nabla _h R_{0cb0} - X_{bh}\nabla _l R^{\_ \_ b}_{lc\_0} \\
& + \omega ^{a}_{bh} \nabla _c R^{\_ \_ b}_{0a\_ 0} - \omega ^{a} _{bh} \nabla _a R^{\_ \_ b}_{0c\_0}
-X_{bc} \nabla ^l R^{\_ \_ b}_{hl\_0} + \omega ^{a} _{bc} \nabla _a R^{\_ \_ b}_{0h\_0} \\
& - \omega ^{a} _{bc} \nabla _h R^{\_ \_ b}_{0a\_0} + \omega ^{a} _{bb} \nabla _c R_{0ha0} - \omega ^{a} _{bb} \nabla _h R_{0ca0} + X_{ba} \nabla _c R^{\_ \_ b}_{0h\_a}  \\
&- X_{ba}\nabla _h R^{\_ \_ b}_{0c\_ a} + (R^{b}_{\_ 0h0}+X^{k} _{b}X_{kh})R^{\_ \_ b}_{0c\_ 0}-(R_{0bah}-X_{b}^k\omega ^{a}_{kh})R^{\_ \_ b}_{ac\_ 0} \\
& + (R_{b0c0}+X_{bk}X_{kc})R^{\_ \_ b}_{h0\_0}
(R_{0bac}-X_{bk}\omega ^{a}_{kc})R^{\_ \_ b}_{ha\_0}   \\
& -(R^{\_ b}_{0\_ab}-X_{bk}\omega ^{a}_{kb})R_{hca0}+(R_{b0a0}+X_{bk}X_{ka})R^{\_ \_ b}_{hc\_a}  \\
= & - \partial _c \nabla _b R^{\_ \_ b}_{0h\_0} + \partial _h \nabla ^b R_{0cb0} -X_{b}^a\nabla _a R^{\_ \_ b}_{hc\_ 0}  - X_{bh}\nabla ^l R^{\_ \_ b}_{lc\_ 0} \\
& + \omega ^{a}_{bh} \nabla _c R^{\_ \_ b}_{0a\_0} - \omega ^{a} _{bh} \nabla _a R^{\_ \_ b}_{0c\_0}
 -X_{bc} \nabla _l R^{\_ \_ b}_{hl\_0} + \omega ^{a} _{bc} \nabla _a R^{\_ \_ b} _{0h\_0} \\
 & - \omega ^{a} _{bc} \nabla _h R^{\_ \_ b}_{0a\_0} + \omega ^{a} _{bb} \nabla _c R_{0ha0} - \omega ^{a} _{bb} \nabla _h R_{0ca0} + X_{ba} \nabla _c R^{\_ \_ b}_{0hb\_a}  \\
 & - X_{ba}\nabla _h R^{\_ \_ b}_{0c\_a} - \partial _b \nabla _c R^{\_ \_ b}_{0h\_0} + \partial ^b \nabla _h R_{0cb0} + \partial _c \nabla ^b R_{0hb0} \\
 &- \partial _h \nabla ^b R_{0cb0}
\end{split}
\label{sep}
\end{align}

At this point we note that
\begin{gather}
\partial _c \nabla ^b R_{0hb0} - \partial ^b \nabla _c R_{0hb0} = v_{cb} ^a \partial _a R_{0hb0} + 2 f_{cb} \nabla ^l R_{lhb0}+ (X_{cb}-X_{bc})\nabla ^l R^{\_ \_ b}_{lh\_0}
\nonumber \\
+(\omega ^{a} _{cb}-\omega ^{a} _{bc})\partial _{a} R^{\_ \_ b}_{0h\_0}  - \partial _c (X_{ba} R^{\_ \_ b}_{a h\_0} + \omega ^{a} _{bh} R^{\_ \_ b}_{0a \_0} + \omega ^{a} _{bb} R_{0ha 0} + X_{ba} R^{\_ \_ b}_{0h\_a})\nonumber \\  + \partial ^b (X_{ca} R_{a hb0} + \omega ^{a} _{ch} R_{0a b0} + \omega ^{a} _{cb}R_{0ha 0} + X_{c}^a R_{0hba})
\nonumber
\end{gather}
Noticing that $\partial _{b} \omega ^{\al} _{ck} - \partial _{c} \omega ^{\al} _{bk}  = R^{\al} _{bck} + \omega ^{\alpha} _{\rho k} (\omega ^{\rho} _{bc} -\omega ^{\rho} _{cb}) -\omega ^{\al} _{b\rho}\omega ^{\rho} _{ck} + \omega ^{\al} _{c\rho}\omega ^{\rho} _{bk}$ up to to the zeroeth order by assumption, we have
\begin{align}
\begin{split}
& \partial _c \nabla ^b R_{0hb0} - \partial ^b \nabla _c R_{0hb0}
= \\
&
v_{cb} ^a \partial _a R^{\_ \_ b}_{0h\_0} + 2 f_{cb} \nabla ^l R^{\_ \_ b}_{lh\_0}+(X_{cb}-X_{bc})\nabla ^l R^{\_ \_ b}_{lh\_ 0}
\\
& +(\omega ^{a} _{cb}-\omega ^{a} _{bc})\nabla _{a} R^{\_ \_ b}_{0h\_0}
-  X_{ba} \partial _c R^{\_ \_ b}_{a h\_0} - \omega ^{a} _{bh} \partial _c R^{\_ \_ b}_{0a \_0} \\
& - \omega ^{a} _{bb} \partial _c R_{0ha 0} - X_{ba} \partial _c R^{\_ \_ b}_{0h\_a} +  X_{ca} \partial ^b R^{a}_{\_ hb0}  \\
 &+ \omega ^{a} _{ch} \partial ^b R_{0a b0} + \omega ^{a} _{cb}\partial ^bR_{0ha 0}+ X_{ca} \partial ^bR_{0hba} \\
 &+(R^{\_ \_ b}_{ah\_0}+R^{\_ \_ b}_{0h\_a})(-R_{0bca}-X_{bp}\omega ^p _{ca}+X_{cp}\omega ^p _{ba}) \\
  &+R^{\_ \_ b}_{0a\_0}(R_{abch}-\omega ^a _{bp} \omega ^p _{ch} + \omega ^a _{cp} \omega ^p _{bh} - X_{ba}X_{ch}+X_{ca}X_{bh})  \\
   &+ R_{0ha0} (R^{\_ \_ \_ b}_{abc\_} - \omega ^a _{bp} \omega ^ p _{cb} - \omega ^a _{cp} \omega ^p _{bb} - X_{ba} X_{c}^b - X_{ca} X_{b}^b)
\end{split}
\label {simpli}
\end{align}
\indent Inserting (\ref{simpli}) and another equation (the same one as (\ref{simpli}) with $h$ and $c$ interchanged and the sign flipped) into equation (\ref{sep}) we see that the zeroeth and the first order terms cancel assuming all the identities hold to order zero.

The other evolution equations are similar. We write only the highest order terms here. The lower order ones (indicated by $L_i$) vanish if we assume (as before) that the identities hold to order zero.
\begin{align}
\begin{split}
\nabla _0 (\nabla ^b R_{0h0b}) =& \, -\partial ^l \nabla ^b R_{lhb0} + L_1  \\
\nabla _0 (\nabla ^b R_{hcbk}-\nabla _0 R_{hc0k})
= &\, \partial _c (\nabla ^l R_{lh0k} + \nabla ^l R_{h0lk})
\\
& -\partial _h (\nabla ^l R_{lc0k}-\nabla _h \nabla ^l R_{c0lk})  + L_2  \\
\nabla _0 (\nabla ^l R_{lh0k}+\nabla ^l R_{h0lk}) = & \,  \partial _{h} (\nabla _{\langle i,}  R_{0c j,k\rangle}) - \partial _{c} (\nabla _{\langle i,}  R_{0h j,k\rangle}) +  L_3   \\
\nabla _0 (\nabla_{\langle i,} R_{0h j,k\rangle}) = & \, \partial^{l} (\nabla _{\langle i,}R_{lh j,k\rangle}) + L_4  \\
\nabla _0 (\nabla _{\langle 0,}R_{hc k,l\rangle}) = & \,  - \partial _c (\nabla _{\langle 0,}R_{0h k,l\rangle}) + \partial _h (\nabla _{\langle 0,}R_{0c k,l\rangle}) + L_5  \\
\nabla _0 (\nabla _{\langle 0,}R_{0h k,l\rangle}) = & \,  -\nabla ^m (\nabla _{\langle 0,}R_{hm k,l\rangle}) +L_6
\end{split}
\label{extra}
\end{align}
 The system (\ref{extra}) (along with equation (\ref{sep}))
 is easily verified to be FOSH with zero as the unique solution if zero initially. \\
\indent We note that $\nabla _h R_{h0\lambda \mu}$ and $\nabla _{\langle i,}R_{j,k\rangle\lambda \mu}$ evolve according to
\begin{align}
\begin{split}
\nabla _0 (\nabla ^h R_{h0\lambda \mu}) &= \, L_7  \\
\nabla _0 (\nabla _{\langle i,}R_{j,k\rangle\lambda \mu}) &= \, L_8
\end{split}
\label{diff}
\end{align}
\indent Finally, we calculate the evolution of $B_{khj}=R_{kh0j}+\partial _k X_{hj} -\partial _h X_{kj} - X_{pj} (\omega ^p _{kh} - \omega ^p _{hk})+X_{kp} \omega ^p _{hj} - X_{hp} \omega ^p _{kj}$ and $W_{hkj} = R_{hkij}-\tilde{R}_{hkij} - X_{ki}X_{jh} + X_{jk}X_{hi}$ (i.e. the definitions of the components of the Riemann tensor)
\begin{align}
\begin{split}
\partial _0 B_{khj} &=\, -X_{aj} B^{\_ \_ a}_{kh\_} + X_h ^p B_{pkj} - X_k ^p B_{phj} + X_{pj} R^{\_ \_ \_ p}_{\langle 0,k,h\rangle\_}  \\
\partial _0 W_{hkj} &= \, -\omega ^i _{lj} R^{\_ \_ \_ p}_{\langle 0,h,k\rangle\_} + \omega ^i _{lj} B^{\_ \_ l}_{hk\_} + X_{h}^{l} W_{klij} - X_{k}^l W_{hlij}
\end{split}
\label{def}
\end{align}
The system (\ref{def}) is FOSH having zero as its solution. If $B$ and $W$ are zero then $R_{\al \be \mu \nu}$ is the Riemann tensor of the metric. \\

\noindent \emph{Proof of theorem \ref{vac}}: Equations
(\ref{evvf}), (\ref{bian1}), (\ref{sep}), (\ref{extra}), (\ref{diff}), (\ref{def}) form a FOSH. Using the uniqueness theory for the same, we conclude that zero is the unique solution (zero in $H^{s-2}$ is the same as zero throughout because of Sobolev embedding) if the variables are zero initially. A calculation shows that  they are zero initially.
Such a calculation is quite long and will not be presented here, but it is done
in essentially the same fashion as in \cite{Dis}. This also holds for the system in lemma (\ref{separation}). This proves that we have a solution to the vacuum Einstein equations satisfying all the conditions required by theorem \ref{vac}. \qed

\subsection{Perfect fluids}
 Just as before, we write equations for the preservation of the gauge. Indeed, we show that the Levi-Civita connection corresponding to the orthonormal frame defined by $\partial _{\al}$ has components $\omega ^{k} _{ij}=\omega ^{k}_{ij}$, $\omega ^{j}_{0i}=0$, $\omega ^{i} _{00} = Y_i$, and $\omega ^{0}_{ij}=X_{ij}$.

Calculations similar to the ones in lemma (\ref{initsymm}) show that the same lemma holds for the Einstein-Euler system as well. We assume this implicitly in what follows. We define $S_{\al \be} = \mathrm{Ricc}_{\al \be} - \rho _{\al \be}$ so that the Einstein equations are $S_{\al \be}=0$.

Explicit computation shows that $[\partial _0, \partial _i] = X_i ^a \partial _a + Y_i \partial _0$, and
$[\partial _i, \partial _j] = (2f_{ij} + X_{[i,j]})\partial _0 + (v^p _{ij}  + \omega ^p _{[i,j]})\partial _p$ where $f=0=v$ initially. If we prove that $f=0=v$ is preserved, then indeed the components of the Levi-Civita connection are as described above. Computations similar to the ones leading to the system (\ref{evvf}) prove that
\begin{align}
\begin{split}
\partial _0 v_{ij} ^p = & \, R^{\_ \_ \_ p}_{0[i,j]\_}+\omega ^p _{k[j,} X_{i]k} + Y^p X_{[i,j]}
\\
& + (v_{ij}^a + \omega _{[i,j]} ^a) X_{a}^p-(v_{a[j,} ^p +\omega ^p_{[a,[j,]})X_{i]a} - \partial _{[i,} X_{j]}^p  \\
2\partial _0 f_{ij} &=\, -R_{0[i,j]0}+2X_{[i,j]}\frac{X_{k}^k}{\mu ^{\prime}} - v_{ij} ^a Y_a  - 2f_{a[j} X _{i]a}-2\widetilde{\nabla} _{[j,} Y_{i]}
\end{split}
\label{evvftwo}
\end{align}
If indeed $v=f=0$, and the first Bianchi identity holds, then both the equations in system (\ref{evvftwo}) are satisfied provided $\mathcal{F}_{ij}=\widetilde{\nabla} _{[i,} Y_{j]}+X_{[i,j]}\frac{X_{kk}}{\mu ^{\prime}}=0$. Let $P_i = Y_i + \partial _i F$. We record the following calculations for further use
\begin{align}
\begin{split}
\cF_{ij} &= \, \widetilde{\nabla} _{[i,} Y_{j]}+  X_{[i,j]} \partial _0 F  \\
&=\, \partial _{[i,}P_{j]}-2f_{ij}\partial _0 F - v_{ij} ^p \partial _p F - \omega ^a _{[i,j]}P_a  \\
\partial _{[k,} \cF _{i]j} &= \, \tilde{L}_1 + \partial _{j} \cF _{hk}
\end{split}
\label{lowerorder}
\end{align}
where $\tilde{L}_i$ denote lower order terms as before. They vanish when all the identities are satisfied. The evolution of $\cF$ is given by
\begin{align}
\begin{split}
\mu ^{\prime} \partial _0 \cF _{ij} &=\, \mathrm{lower \ order} + \partial _ {[i}(\partial _0 P_{,j]})  \\
&=\, \mathrm{lower \ order} + \partial_{[i,}\partial ^l X_{j]l} - \partial _{[i,}\partial _{j]}X_{l}^l  \\
&=\, \mathrm{lower \ order} + \partial_{l} \partial _{[i,} X_{j]l}  \\\
&=\, \tilde{L}_2 + \partial ^l B_{ijl}
\end{split}
\label{thefone}
\end{align}
We now write the evolution equation of $R_{\langle \al , \be, \mu \rangle \nu}$ as before:
\begin{align}
\begin{split}
\nabla _0 R_{\langle h,c,d \rangle\al} &=\,  -\nabla _h R_{\langle c,0,d\rangle\al}
 - \nabla _d R_{\langle h,0,c\rangle\al} - \nabla _c R_{\langle 0,h,d\rangle\al}
 \\
&
+ \nabla _{\langle h,}R_{d,c\rangle 0\al}  \\
\nabla _0 R_{\langle 0,h,c\rangle\al} &= \, \nabla ^l R_{\langle l,h,c\rangle}-\nabla ^{\langle l,} R_{h,c\rangle l\al} + \nabla_c S_{h\al}-\nabla_h S_{c\al}
\end{split}
\label{bianevol}
\end{align}
We wish to make sure that the Euler equation $Y_i = -\partial _i F$ is satisfied. The evolution of $P_i$ is computed to be
\begin{align}
\begin{split}
\partial _0 P_i &= \, \frac{1}{\mu ^{\prime}}\left [ \widetilde{\nabla} ^j X_{ij} - \widetilde{\nabla} _i X_{k}^k - Y^k (X_{ki} -X_{ik}) \right ] -X_{i}^l P_l
 \\
&= \, \frac{1}{\mu ^{\prime}} \left [ S_{i0}+ (-\mathrm{Ricc}_{i0}+\widetilde{\nabla} ^j X_{ij} - \widetilde{\nabla} _i X_{k}^k - Y^k (X_{ki} -X_{ik})) \right ]
\\
& - X_{i}^l P_l
\end{split}
\label{fluidone}
\end{align}
Now, we calculate the evolution of the (quasi-)constraints (remembering that $R_{h0i0}= \widehat{R_{h0i0}}-\cF _{hi}$ where $\widehat{R_{h0i0}}$ is the ``true" Riemann tensor). Let $B_{khj}=R_{kh0j}+\partial _k X_{hj} -\partial _h X_{kj} - X_{pj} (\omega ^p _{kh} - \omega ^p _{hk})+X_{kp} \omega ^p _{hj} - X_{hp} \omega ^p _{kj}-Y_j(X_{kh}-X_{hk})$ and $W_{hkij} = R_{hkij}-\tilde{R}_{hkij} - X_{ki}X_{jh} + X_{jk}X_{hi}$.
\begin{align}
\begin{split}
 \partial _0 B_{khj}   = &\,  -X_{lj}^{\_ \_ l}B_{kh\_} + X_{k}^l B_{hlj} + X_{h}^l B_{lkj} + X_{pj} R^{\_ \_ \_ p}_{\langle 0,k,h\rangle\_}
 \\
 & -Y^p W_{hkpj}- Y_j R_{0[k,h]0}
  +2f_{kh}\partial _0 Y_j
  \\
  & + v^p _{kh} \partial _p Y_j +Y_{[k,} \cF _{h]j} -Y_j \cF _{[k,h]}+\widetilde{\nabla}_{[k,} \cF _{h]j}  \\
   = & \,  \tilde{L}_3 + \partial _{j} \cF _{hk}  \\
 \partial _0 W_{hkij}  = & \, \tilde{L}_4  \\
 \partial _0 \nabla _{\langle i,} R_{j,k\rangle \mu \nu} = & \, \tilde{L}_5  \\
 \partial _0 (\nabla _h R_{h0\lambda \mu}  -
 \,\,\,
 &
\nabla _{\lambda} \rho _{\mu 0} +
 \nabla _{\mu}  \rho _{\lambda 0})  = \, \tilde{L}_6
\end{split}
\label{constr2}
\end{align}
\indent Finally, we compute the evolution of $S_{\al \be}$
\begin{align}
\begin{split}
\nabla _0 S_{ab} &= \, \nabla _a S_{0b} + \nabla ^h R_{\langle 0,a,h\rangle b}-(\nabla ^h R_{h0ab}-\nabla _a \rho_{b0}+\nabla _b \rho _{a0})  \\
\nabla _0 S_{0a} &= \, \nabla ^k S_{ka} +2(\mu+p)P_b  \\
\nabla _0 S_{a0} &= \, \nabla_a S_{00} + \nabla ^h (R_{0h0b}-R_{0b0h})  \\
\nabla _0 S_{00} &= \, \nabla ^k S_{k0}
\end{split}
\label{Einstevolve}
\end{align}
Notice that (\ref{evvftwo}), (\ref{thefone}), (\ref{bianevol}), (\ref{fluidone}), (\ref{constr2}), and (\ref{Einstevolve}) form a FOSH. They have zero as their unique solution if the Bianchi identities, the constraints, and the Euler equations hold.

\begin{remark} Notice that the preservation of gauge and the satisfaction of Einstein's equations maybe proven in the vacuum case in a manner similar to that of the fluid case. Actually, the vacuum case maybe deduced from the fluid case by
keeping track of Newton's constant $G$ in the equations and setting $G=0$. However, we chose to do it otherwise in order to separate the preservation of the gauge and constraints from the Einstein equations themselves. This maybe useful in other contexts.
\end{remark}

\noindent \emph{Proof of theorem \ref{EE}}: By arguments similar to those used in the proof of theorem \ref{vac}, we obtain a solution of the Einstein-Euler system. This satisfies almost all the conditions required by theorem \ref{EE} except ostensibly, the regularity of $g$, because of the regularity of $b_i$ that is lower than desired. This however, is an artifact of our chosen coordinate system  (which depends on $v$ which in turn has lower regularity than $h_0$). To see this, first we note that the $g$ we obtained is in $C^{0}([0,T],H^{s}(U))\cap C^{1}([0,T],H^{s-1}(U)) \cap C^{2}([0,T],H^{s-2}(U))$,
with $U$ some local chart (see remark \ref{remark_charts}). However, its restriction to $\Sigma \times \{0\}$ is in $H^{s+1}$. Using the ADM decompositon of lapse and shift, one may choose the initial lapse to be $1$, the initial shift vector to be $0$ and their time derivative appropriately so as to satisfy the wave gauge condition initially. The energy-momentum tensor will be in $C^{0}([0,T],H^{s}(U))\cap C^{1}([0,T],H^{s-1}(U))$. The Einstein equations form a quasilinear hyperbolic system in the wave gauge. Hence we get a metric solving the full Einstein equations with the correct regularity. This coincides our original solution due to local geometric uniqueness \cite{CBbook}. \qed


\end{document}